\documentclass[12pt]{amsart}
\usepackage{amsmath, amsthm, euscript}

\def\nxn{n{\times}n}
\def\T{\overline{T}}
\def\Utilde{\widetilde{U}}
\def\RelMats{\mathit{RelMatrices}}
\def\RelEnts{\mathit{RelEntries}}
\def\reply{\mathit{reply}}
\def\set{\mathit{set}}

\def\Mon{\mathit{Monomials}}
\def\Monn{\mathit{Monomials}_n}
\def\Monnbar{\overline{\mathit{Monomials}_n}}
\def\Monbar{\overline{\mathit{Monomials}}}
\def\Charn{\mathit{CharCoefs}_n}
\def\Char{\mathit{CharCoefs}}
\def\Charnbar{\overline{\mathit{CharCoefs}_n}}
\def\Charbar{\overline{\mathit{CharCoefs}}}
\def\RelId{\mathit{RelIdeal}}
\def\trace{\operatorname{trace}}
\def\det{\operatorname{det}}
\def\nxn{n{\times}n}
\def\mxm{m{\times}m}
\def\bw{\mathbf{w}}
\def\bx{\mathbf{x}}
\def\by{\mathbf{y}}
\def\bxbar{\overline{\mathbf{x}}}

\def\X{\overline{X}}

\def\A{\overline{A}}
\def\B{\overline{B}}
\def\M{\overline{M}}
\def\Tn{\overline{T}_n}

\def\Sn{\overline{S}_n}
\def\sltwo{\mathfrak{sl}_2}
\def\Q{\mathbb{Q}}
\def\F{\mathbb{F}}

\def\PSLtwoZ{\mathit{PSL}_2(\mathbb{Z})}

\swapnumbers
\theoremstyle{plain}
\newtheorem{thm}{Theorem}[section]
\newtheorem{lem}[thm]{Lemma}

\theoremstyle{definition}

\newtheorem*{ack}{Acknowledgement}

\theoremstyle{remark}
\newtheorem{note}[thm]{}

\def\n{\noindent}

\begin{document}

\title[Detecting Infinitely Many Semisimple
Representations]{Detecting Infinitely Many Semisimple Representations
  in a Fixed Finite Dimension}

\author{E. S. Letzter}

\address{Department of Mathematics\\
        Temple University\\
        Philadelphia, PA 19122-6094}
      

      \email{letzter@temple.edu }

      \thanks{Research supported in part by grants from the National Security
        Agency.}



\begin{abstract} Let $n$ be a positive integer, and let $k$ be a field
  (of arbitrary characteristic) accessible to symbolic computation. We
  describe an algorithmic test for determining whether or not a
  finitely presented $k$-algebra $R$ has infinitely many equivalence
  classes of semisimple representations $R \rightarrow M_n(k')$, where
  $k'$ is the algebraic closure of $k$. The test reduces the problem
  to computational commutative algebra over $k$, via famous results of
  Artin, Procesi, and Shirshov. The test is illustrated by explicit
  examples, with $n = 3$.
\end{abstract}

\maketitle

\section{Introduction}

Among the most fundamental tasks, when studying a given finitely
presented algebra over a field $k$, is the parametrization of the
irreducible finite dimensional representations. Typically, such
parametrizations depend on whether or not the finite dimensional
irreducible representations, partitioned according to their
dimensions, occur in finite or infinite families. The focus in this
paper is on general algorithmic approaches to this latter issue.

\begin{note} Before describing the work of this paper, we review some
  of the background and context. To start, assume that $k$ is a field
  accessible to symbolic computation and that $n$ is a positive
  integer. Set
\[R := k\{ X_1,\ldots, X_s \}\big/\langle f_1,\ldots, f_t \rangle, \]
the free associative $k$-algebra in $X_1,\ldots,X_s$ modulo the ideal
generated by $f_1,\ldots,f_t$. By an \emph{$n$-dimensional
  representation of $R$\/} we will always mean a unital $k$-algebra
homomorphism
\[ R \rightarrow M_n(k') , \]
where $k'$ is the algebraic closure of $k$. Irreducibility,
semisimplicity, and equivalence of representations are defined over
$k'$; see (\ref{IIr2}).

Since the presence of finite dimensional representations (without
restriction on the dimension) over $k$ is a Markov property
\cite{Mar},\cite{Rab} the question of existence of finite dimensional
representations of $R$ is algorithmically undecidable in general
\cite{Bok}. Consequently, the task of algorithmically studying the
finite dimensional representation theory of $R$, without bounding the
dimensions involved, appears to be hopeless. However, if we fix $n$
and restrict our attention to dimensions $\leq n$, then the
situation improves considerably (at least in principle). To start, the
existence of representations in dimensions bounded by $n$ is
determined by finitely many (commutative) polynomial equations, and so
can be approached (again in principle) using computational commutative
algebra. In \cite{Let3} tests determining the existence of irreducible
$n$-dimensional representations were described, and in \cite{Let1}
tests determining the existence of non-semisimple
at-most-$n$-dimensional representations of $R$ were described.
\end{note} 

\begin{note} Our purpose in this paper is to describe a test,
  involving computational commutative algebra over $k$, for
  determining whether or not $R$ has infinitely many equivalence
  classes of $n$-dimensional semisimple representations. The test is
  developed in \S \ref{II} and presented in \S \ref{III}.
\end{note}

\begin{note}
  In \S \ref{IV}, explicit illustrative examples are given, with $n=3$
  and $s = 2$. The calculations in these examples were performed with
  Macaulay2 \cite{Mac} on small computers ($\leq$ 8 GB RAM).
\end{note}

\begin{note} Our approach can be sketched as follows. To start, by
  famous results of Artin \cite{Art} and Procesi \cite{Pro}, the
  equivalence classes of $n$-dimensional representations of $R$
  correspond exactly with the maximal ideals of the (commutative)
  trace ring of $R$, taken over $k'$; see (\ref{IIr8}). In particular,
  $R$ has finitely many equivalence classes of $n$-dimensional
  semisimple representations if and only if this trace ring is finite
  dimensional over $k'$. Next, using Shirshov's Theorem \cite{Shi},
  the trace ring's being finite or infinite dimensional depends only
  on a suitably truncated trace ring. (We use Belov's refinement of
  Shirshov's Theorem \cite{Bel}.) Finally, the finite-versus-infinite
  dimensionality of the truncated trace ring can be algorithmically
  determined using a variant of the subring membership test, working
  over $k$.
\end{note}

\begin{note} An algorithm, in characteristic zero, for determining
  whether or not there exist infinitely many equivalence classes of
  $n$-dimensional irreducible representations was outlined in
  \cite{Let2}. In part, the present paper provides a generalization
  and simplification of this previous work. Moreover, the methods in
  \cite{Let2} can be combined with the approach of the present paper
  to formulate a test for determining whether or not there exist
  infinitely many equivalence classes of irreducible $n$-dimensional
  representations of $R$, in arbitrary characteristic. However, the
  methods in \cite{Let2} appear to be considerably more costly than
  the test described in the present paper. A detailed analysis, with
  examples, of tests for detecting infinitely many equivalence classes
  of irreducible $n$-dimensional representations is left for future
  work.
\end{note}

\begin{note} A warning: Since we are assuming that representations are
  unital maps, it is possible for $R$ to have only finitely many
  (possibly zero) equivalence classes of $n$-dimensional semisimple
  representations but at the same time for $R$ to have infinitely many
  equivalence classes of irreducible representations in some dimension
  less than $n$. Of course, if $R$ has a $1$-dimensional
  representation, then finiteness of the number of equivalence classes
  of $n$-dimensional semisimple representations ensures the finiteness
  of the number of equivalence classes of irreducible representations
  in dimensions $\leq n$.
\end{note}

\begin{ack} The author is grateful to the referee for suggestions
on improving the exposition.
\end{ack}

\section{Setup and Proof of Test} \label{II}

In this section we develop and prove our test to determine whether a
finitely presented algebra over a field has infinitely many distinct
equivalence classes of $n$-dimensional semisimple representations, for
a fixed $n$. The algorithm is presented in \S \ref{III}.

\begin{note} \label{IIr1} Assume that $n$ is a positive integer, that $k$ is
a field, and that $k'$ is the algebraic closure of $k$. Let $k\{ X_1,\ldots,
X_s \}$ denote the free associative $k$-algebra in the noncommuting
indeterminates $X_1,\ldots,X_s$.

Set 
\[R := k\{ X_1,\ldots, X_s \}\big/\langle f_1,\ldots, f_t \rangle,\]
for some fixed choice of $f_1,\ldots , f_t$ in $k\{X_1,\ldots,X_s
\}$. We will use $\X_1,\ldots,\X_s$ to denote the respective images in $R$ of
$X_1,\ldots,X_s$.
\end{note}

\begin{note} \label{IIr2} Let $L$ be any subfield of $k'$, let $m$ be a positive
integer, and let $\Lambda$ be an $L$-algebra. 

(i) By an \emph{$m$-dimensional representation} of $\Lambda$ we will always
mean a unital $L$-algebra homomorphism from $\Lambda$ into the $L$-algebra
$M_m(k')$ of $\mxm$ matrices over $k'$.

(ii) Representations $\rho, \rho' \colon \Lambda \rightarrow M_m(k')$ are
\emph{equivalent} if there exists a matrix $Q \in GL_m(k')$ such that
\[\rho'(a) = Q\rho(a) Q^{-1},\]
for all $a \in \Lambda$. 

(iii) An $m$-dimensional representation $\rho \colon \Lambda \rightarrow
M_m(k')$ is \emph{irreducible} (cf.~\cite[\S 9]{Art}) provided
$k'\rho(\Lambda) = M_m(k')$.

(iv) An $n$-dimensional representation $\rho \colon \Lambda \rightarrow
M_n(k')$ is \emph{semisimple} provided $\rho$ is equivalent to a
representation of the form
\[a \; \longmapsto \; \bmatrix \rho_1(a) \\ & \rho_2(a) \\ & & \ddots \\ & & &
\rho_r(a) \endbmatrix,\]
for suitable choices of positive integers $m_1,\ldots,m_r$, suitable
choices of irreducible $m_i$-dimensional representations $\rho_i$ of
$\Lambda$, and all $a \in \Lambda$.
\end{note}

\begin{note} \label{IIr3} (i) Set
\[B' := k' [ x_{ij}(\ell) : 1 \leq i,j \leq n, \, 1 \leq \ell \leq s ],\]
and
\[B := k [ x_{ij}(\ell) : 1 \leq i,j \leq n, \, 1 \leq \ell
\leq s ] \; \subseteq \; B',\]
where the $x_{ij}(\ell)$ are commuting indeterminates.

For $1 \leq \ell \leq s$, let $\bx _\ell$ denote the $\nxn$ generic
matrix $(x_{ij}(\ell))$, in $M_n(B)$. For $g = g(X_1,\ldots,X_s) \in
k\{ X_1, \ldots , X_s \}$, let $g(\bx) = g(\bx_1,\ldots,\bx_s)$ denote
the image of $g$, in $M_n(B) \subseteq M_n(B')$, under the canonical
map
\[ k\{ X_1,\ldots , X_s \} \xrightarrow{\; X_\ell \; \longmapsto \; \bx_\ell
\; } M_n(B) . \]
Identify $B'$ with the center of $M_n(B')$,
and identify $B$ with the center of $M_n(B)$.

(ii) Set 
\[\RelMats \; := \; \{ f_1(\bx), \ldots , f_t (\bx)
\}.\]
Let $\RelId(M_n(B'))$ be the ideal of $M_n(B')$ generated by
$\RelMats$, and let $\RelId(M_n(B))$ be the ideal of $M_n(B)$ generated
by $\RelMats$. Note that
\[\RelId(M_n(B)) \; = \; \RelId(M_n(B')) \cap M_n(B) .\]

(iii) Let $\RelEnts$ denote the set of entries of the matrices in
$\RelMats$. Let $\RelId(B)$ denote the ideal of $B$ generated by
$\RelEnts$, and let $\RelId(B')$ denote the ideal of $B'$ generated by
$\RelEnts$. Then
\[\RelId(B) \; = \; \RelId(B') \cap B.\]
Also,
\[\RelId(B') \; = \; \RelId(M_n(B')) \cap B', \quad \text{and} \quad
\RelId(B) \; = \; \RelId(M_n(B)) \cap B.\]

(iv) Set
\[A' \; := \; k'\{ \bx_1, \ldots , \bx_s \}, \]
the $k'$-subalgebra of $M_n(B')$ generated by the generic matrices $\bx_1,
\ldots , \bx_s$, and set
\[A \; := \; k\{ \bx_1, \ldots , \bx_s \} \; \subseteq \; A'.\]
Set
\[\RelId (A') \; := \; \RelId (M_n(B')) \cap A',\]
and
\[ \RelId (A) \; := \; \RelId (M_n(B)) \cap A.\]
\end{note}

\begin{note} \label{IIr4} Let 
\[\A  :=  A/\RelId(A), \quad \B  :=  B/\RelId(B), \quad \M  := 
M_n(B)\big/\RelId \big(M_n(B)\big)\]
and
\[\A'  :=  A'/\RelId(A'), \; \B'  :=  B'/\RelId(B'), \; \M'  :=
 M_n(B')\big/\RelId \big(M_n(B')\big).\]
 Identify $\A$, $\B$, and $\M$ with their respective natural images in
 $\A'$, $\B'$, and $\M'$. Since
\[\RelId(M_n(B')) \; = \; M_n(\RelId(B')),\]
we see that there is a natural isomorphism
\[\M' \; = \; M_n(B')/M_n(\RelId(B')) \; \cong \; M_n(\B') .\]
We use this isomorphism to identify $\M'$ with $M_n(\B')$.
Furthermore, $\B'$ is isomorphic to the image of $B'$ in $\M'$; we
identify $\B'$ with that image. In particular, $\B' = Z(\M')$, the
center of $\M'$. We similarly identify $\B$ with $Z(\M)$.

 Let $\bxbar_1,\ldots,\bxbar_s$ denote, respectively, the images of
 $\bx_1,\ldots,\bx_s$ in $\M \subseteq \M'$. So
 $\bxbar_1,\ldots,\bxbar_s$ generate $\A'$ as a $k'$-algebra and
 generate $\A$ as a $k$-algebra.
\end{note}

\begin{note} \label{IIr5} We have a canonical $k$-algebra homomorphism
\[\pi' : R \; \xrightarrow{\quad \X_\ell \; \longmapsto \; \bxbar_\ell \quad}
\; \M \; \xrightarrow{\quad \text{inclusion} \quad} \; \M'.\]
Note that $\pi'(R) = \A$ and that $\A'$ is generated as a $k'$-algebra by
$\pi'(R)$. 
\end{note}

\begin{note} \label{IIr6} We will say that a $k'$-algebra homomorphism $\alpha
\colon \M' \rightarrow M_n(k')$ is \emph{matrix-unital} provided $\alpha$
restricts to the identity on $M_n(k') \subseteq \M'$. More generally, if $L$
is a subfield of $k'$ and $\Lambda$ is an $L$-subalgebra of $\M'$, we will say
that an $n$-dimensional representation $\rho \colon \Lambda \rightarrow
M_n(k')$ is \emph{matrix-unital (with respect to $\M'$)} when $\rho$ is the
restriction of a matrix-unital map $\M' \rightarrow M_n(k')$. In other words,
the matrix-unital $n$-dimensional representations of $\Lambda$ all have the
form
\[\Lambda \; \xrightarrow{\quad \text{inclusion} \quad} \; \M' \;
\longrightarrow \; M_n(k'),\]
where the right-hand map is matrix-unital. 
\end{note}

\begin{note} \label{IIr7} Every $n$-dimensional representation $\rho:R
\rightarrow M_n(k')$ can be written in the form
\[R \; \xrightarrow{\quad \pi' \quad} \; \M' \; \xrightarrow{\; \rho^{\M'} \;}
\; M_n(k'),
\]
for a suitable, unique, matrix-unital representation $\rho^{\M'}:\M'
\rightarrow M_n(k')$. (Conversely, of course, every matrix-unital
$n$-dimensional representation of $\M'$ gives rise to a unique
$n$-dimensional representation of $R$, written in the preceding form.)
Since $\A'$ is the $k'$-subalgebra of $\M'$ generated by $\pi'(R)$, we
see that the assignment
\[ \rho \; \mapsto \; \rho^{\A'} \; := \; \rho^{\M'}\vert _{\A'} \]
induces a one-to-one correspondence between the equivalence classes of
$n$-dimensional representations of $R$ and the equivalence classes of
$n$-dimensional matrix-unital representations of $\A'$. This correspondence
depends only on the original choice of presentation of $R$ given in
(\ref{IIr1}).

Furthermore, an $n$-dimensional representation $\rho:R \rightarrow M_n(k')$
will be semisimple if and only if the corresponding matrix-unital
representation $\rho^{\A'}$ is semisimple, and $\rho$ will be irreducible
if and only if $\rho^{\A'}$ is irreducible.
\end{note}

\begin{note} \label{IIr8} (i) Let $\Monbar \subseteq \M \subseteq \M'$
  denote the set of monomials (i.e., matrix products) of length
  greater than or equal to $1$, in the $\bxbar_1,\ldots,\bxbar_s$. Let $\Charbar
  \subseteq \B \subseteq \B'$ denote the set of nonscalar coefficients
  of characteristic polynomials of monomials in $\Monbar$. Let $\T$
  denote the $k$-subalgebra of $\B$ generated by $\Charbar$, and let
  $\T'$ denote the $k'$-subalgebra of $\T'$ generated by
  $\Charbar$. It follows from Shirshov's theorem \cite{Shi} (see
  \cite[Proposition 3.1]{Pro}) that $\T'$ is a finitely generated
  $k'$-algebra. 

  In the literature, $\T'$ is referred to as a \emph{trace ring},
  since in characteristic zero it is generated by traces.  More
  precisely, in characteristic zero, Razmyslov proved that $\T'$ is
  generated by the traces of the monomials in $\Monbar$ of length
  $\leq n^2$; this upper bound is the best known \cite{Raz}.

  (ii) Observe that any matrix-unital representation of $\M'$ will map
  coefficients of characteristic polynomials of matrices in $\M'$ to
  coefficients of characteristic polynomials of matrices in $M_n(k')$.

  (iii) Now let $\rho:R \rightarrow M_n(k')$ be an arbitrary
  $n$-dimensional representation, with corresponding matrix-unital map
  $\rho^{\M'}:\M' \rightarrow M_n(k')$. By (ii), the restriction of
  $\rho^{\M'}$ to $\T'$ produces a $k'$-algebra homomorphism
  $\rho^{\T'}:\T' \rightarrow k'$.  Since coefficients of
  characteristic polynomials are invariant under conjugation, we
  further see that $\rho^{\T'}$ depends only on the equivalence class
  of $\rho$. Consequently, the assignment $\rho \mapsto \rho^{\T'}$
  provides a well defined function from the set of equivalence classes
  of $n$-dimensional representations of $R$ to the set of $k'$-algebra
  homomorphisms from $\T'$ onto $k'$. Moreover, the assignment $\rho
  \mapsto \ker\rho^{\T'}$ then provides a well defined function from
  the set of equivalence classes of $n$-dimensional representations of
  $R$ to $\max \T'$.

  (iv) Key to Artin's \cite{Art} and Procesi's \cite{Pro} study of
  finite dimensional representations is their proof that the function
  $\rho \mapsto \ker\rho^{\T'}$ induces a bijection from the set of
  equivalence classes of semisimple $n$-dimensional representations of
  $R$ onto the maximal spectrum of $\T'$.

In particular, there are only finitely many equivalence classes of
$n$-dimensional semisimple representations of $R$ if and only if the finitely
generated commutative $k'$-algebra $\T'$ is finite dimensional over $k'$.
\end{note}

\begin{note} \label{IIr9} 
Let $\Monnbar$ denote the set of monomials in $\Monbar$ of length less than or
equal to $n$. (The use of the number $n$ will be explained below.) Let
$\Charnbar$ denote the set of nonscalar coefficients of characteristic
polynomials of monomials in $\Monnbar$. Let $\Tn$ denote the
$k$-subalgebra of $\B$ generated by $\Charnbar$, and let $\Tn'$ denote
the $k'$-subalgebra of $\B'$ generated by $\Charnbar$.

Set
\[\Sn' \; = \; \Tn'\{ \bxbar_1,\ldots, \bxbar_s \}.\]
We can conclude as follows that $\Sn'$ is a finitely generated module
over the central subalgebra $\Tn'$: To start, by Belov's refinement
\cite{Bel} of Shirshov's Theorem (see \cite[\S 9.2]{Dre}), it follows
that there exists an integer $h$, depending only on $n$ and $s$,
such that $\A'$ is spanned as a $k'$-vector space by the products
$w_1^{a_1}\cdots w_m^{a_m}$, where $m \leq h$ and where
$w_1,\ldots,w_m$ are monomials of degree $\leq n$. (In Shirshov's
original theorem the degree bound on $w_1,\ldots,w_m$ is $2n-1$.) By the
Cayley-Hamilton Theorem, each $w_i^{a_i}$, for $a_i \geq n$, can be
expressed as a $\Tn'$-linear combination of powers of $w_i$ of degree
less than $n$.  Hence $\Sn'$ is a finitely generated
$\Tn'$-module. The length $n$ used in defining $\Charnbar$ and $\Tn'$
is the minimum length required to apply Belov's result.
\end{note}

The following is now a corollary to (\ref{IIr8}) and (\ref{IIr9})

\begin{lem} \label{IIt1} $R$ has only finitely many  equivalence classes
of semi\-simple $n$-di\-men\-sion\-al representations if and only if $\Tn'$
is a finite dimensional $k'$-algebra.
\end{lem}

\begin{proof} Suppose first that $R$ has only finitely many equivalence
classes of semisimple $n$-dimensional representations. By (\ref{IIr8}iv), $\T'$
is finite dimensional over $k'$, and so $\Tn'$ is finite dimensional over
$k'$.

Next, suppose that $\Tn'$ is finite dimensional over $k'$. By (\ref{IIr9}),
$\Sn'$ is finite dimensional over $k'$, and so $\A' \subseteq \Sn'$ is finite
dimensional over $k'$. Consequently, it follows from (\ref{IIr7}) that $R$ has
only finitely many equivalence classes of $n$-dimensional semisimple
representations.
\end{proof}

\begin{note} \label{IIr10} By (\ref{IIt1}), to algorithmically decide whether or
not $R$ has infinitely many distinct equivalence classes of $n$-dimensional
semisimple representations it remains to find an effective means of
determining when $\Tn'$ is finite dimensional over $k'$. To start, observe that
$\Tn'$ is finite dimensional over $k'$ if and only if each $d \in \Charnbar$ is
algebraic over $k'$, if and only if each $d \in \Charnbar$ is algebraic over
$k$ (since $k'$ is algebraic over $k$).

Now let $\Monn \subseteq M \subseteq M'$ denote the set of monomials (i.e.,
matrix products) of length greater than or equal to $1$, but less than or
equal to $n$, in the $\bx_1,\ldots,\bx_s$. Let $\Charn \subseteq B \subseteq
B'$ denote the set of nonscalar coefficients of characteristic
polynomials of monomials in $\Monn$. It follows from the preceding paragraph
that $R$ has at most finitely many equivalence classes of $n$-dimensional
semisimple representations if and only if each $c \in \Charn$ is algebraic,
modulo $\RelId(B)$, over $k$.  We can use the following variant of the subring
membership test to determine whether a given $c \in \Charn$ is algebraic,
modulo $\RelId(B)$, over $k$.
\end{note}

\begin{note} \label{IIr11} (Cf.~\cite[pp.~269--270]{BecWei}.) Let $U$ denote
the commutative polynomial ring $k[y_1,\ldots,y_m]$, let $a_1,\ldots,a_u \in
U$, and let $I$ be the ideal of $U$ generated by $a_1,\ldots,a_u$. Choose $g
\in U$. Then $g$ is algebraic, modulo $I$, over $k$ if and only if $I \cap
k[g] \ne 0$. Now view $U$ as a subring of $\Utilde = k(v)[y_1,\ldots,y_m]$,
where $v$ is an indeterminate. Then $I \cap k[g] \ne 0$ if and only if $1$ is
contained in the ideal
\[J \; := \; (v-g).\Utilde + I.\Utilde \quad = \quad \langle v-g, a_1,\ldots,a_u \rangle\]
of $\Utilde$.

Assuming that the ideal membership test can be applied to (commutative)
polynomial rings over $k(v)$, we can effectively determine whether or
not $1$ is contained in $J$, and we can thus determine whether or not
$g \in U$ is algebraic, modulo $I$, over $k$.
\end{note}

\begin{note} \label{IIr12} In characteristic $0$, we may replace
  $\Charn$ in (\ref{IIr10}) with the set of traces of elements of
  $\Mon$ of length $\leq n^2$. This replacement is possible because,
  in characteristic zero, $\T'$ is generated by the traces of those
  elements of $\Monbar$ of length $\leq n^2$, following (\ref{IIr8}i).
  Details are left to the interested reader.
\end{note}

\section{Test for Detecting Infinitely Many Semisimple\\
  Representations in a Fixed Finite Dimension}
\label{III}

Retain the notation of the preceding section (although we have
attempted to make the discussion below reasonably self contained). It
immediately follows from (\ref{IIr10}) and (\ref{IIr11}) that the
following procedure will determine whether or not the finitely
presented algebra $R$ has infinitely many distinct equivalence classes
of semisimple $n$-dimensional representations.

\smallskip\n\textbf{Inputs:} A positive integer $n$, a field $k$ (suitably 
accessible to symbolic computations) with algebraic closure $k'$,
a finitely presented $k$-algebra
\[R = k\{ X_1,\ldots, X_s \}\big/\big\langle \, f_1(X_1,\ldots,X_s),\ldots,
f_t(X_1,\ldots,X_s) \, \big\rangle\]

\n\textbf{Output:} YES if there are infinitely many distinct
equivalence classes of semisimple representations $R \rightarrow M_n(k')$; NO
otherwise

\begin{tabbing}
\hspace{3ex}\=\hspace{3ex}\=\hspace{3ex}\=\hspace{3ex}\=\kill

\textbf{begin} \\[4pt]

\> $C := k(v) [ x_{ij}(\ell) : 1 \leq i,j \leq n, 1 \leq \ell \leq s ]$
\\[4pt]

\> $\RelEnts :=$ set of entries of the matrices $f_1(\bx_1,\ldots,\bx_s),
\ldots, f_t(\bx_1,\ldots,\bx_s)$, \\

\> where $\bx_\ell$ denotes the $\nxn$ generic matrix $(x_{ij}(\ell))$, 
for $1 \leq \ell \leq s$ \\[4pt]

\> $\Monn :=$ the set of matrix products of length greater than or equal to $1$,
\\

\>  but less than or equal to $n$, in the $\bx_1,\ldots,\bx_s$
\\[4pt]

\> $\Charn :=$ the set of nonscalar coefficients of characteristic polynomials\\

\> of matrices in $\Monn$ \\[4pt]

\>  $\reply :=$ ``NO''
\\[4pt]

\> $\set := \Charn$
\\[4pt]

\> \> \textbf{while} $\set \ne \emptyset$ \textbf{do}
\\[4pt]

\> \> Choose $c \in \set$ \\[4pt]

\> \>  $J :=$ ideal of $C$ generated by $c-v$ and $\RelEnts$
\\[4pt]

\> \> \> \textbf{If} $1 \notin J$ \; (applying Ideal Membership Test) 
\\[4pt]

\> \> \> \> \textbf{then} $\reply :=$ ``YES'' \textbf{and} $set := \emptyset$ 
\\[4pt]

\> \> \> \> \textbf{else} $\set := set \setminus \{ c \}$
\\[4pt]

\> \> \> \textbf{end}
\\[4pt]

\> \> \textbf{end}
\\[4pt]

\> \textbf{return} $\reply$
\\[4pt]

\textbf{end} 
\end{tabbing}

\begin{note} To give some relative measure of the complexity of the
  above algorithm, note that each generator of the ideal $J$ will be a
  (commutative) polynomial, over $k(v)$, of degree no greater than
  $\max(n^2,e)$, where $e$ is the maximum total degree of the
  $f_1,\ldots,f_t$. In particular, the members of $\Charn$ used in the
  algorithm are polynomials of degree $\leq n^2$. 
 \end{note}

\section{Examples: $3$-dimensional representations of $2$-generator algebras}
\label{IV}

Retain the notation of the previous sections. To illustrate the procedure in
\S \ref{III}, we examine the $3$-dimensional representations of algebras of
the form
\[ R \; = \; k\{ X , Y \} \big/ \langle f_1(X,Y), \ldots , f_t(X,Y)
\rangle. \]
To start, set $\bx = (x_{ij})$ and $\by = (y_{ij})$, in $M_3(C)$, where
\[ C \; = \; k(v)[x_{ij}, y_{ij} : 1 \leq i,j \leq 3]. \]
Then $\RelEnts$, in this situation, is a set of $9$ (commutative)
polynomials, in $9t$ variables, with coefficients in $k$. The maximum
degree appearing is the maximum of the total degrees of the
$f_1,\ldots,f_t$. Also, $\Mon_3$ consists of the $14$ monomials in
$\bx$ and $\by$ of length $1$, $2$, or $3$.

The characteristic polynomial of a $3{\times}3$ matrix $\bw = (w_{ij})$ is
\[ \lambda^3 -\trace(\bw)\lambda + (w_{11}w_{22} + w_{11}w_{33} + w_{22}w_{33}
-w_{12}w_{21}-w_{13}w_{31}-w_{23}w_{32})\lambda - \det(\bw) ,\]
and so $\Char_3$ is a set of $42$ distinct polynomials, in $18$
variables, with coefficients in $k$. The maximum degree appearing is
$9$.

To perform the algorithm in \S \ref{III}, we need to check, for $c \in
\Char_3$, whether $1$ is contained in the ideal of $C$ generated by $(c-v)$
and $\RelEnts$. If for all $c \in \Char_3$ these ideals contain $1$, then
there are only finitely many equivalence classes of semisimple $3$-dimensional
representations $R \rightarrow M_3(k')$. If there is at least one of these
ideals that does not contain $1$, then there are infinitely many equivalence
classes.

\begin{note}
For a concrete example, let
\[ R \; = \; k\{ X, Y \} \big/ \langle X^2 - 1, \; Y^3 - 1 \rangle ,\]
the group algebra, over $k$, of $\PSLtwoZ$. Using Macaulay2
\cite{Mac}, for $k = \Q$, $\F_2$, $\F_3$, $\F_5$, and $\F_7$, we found
that $1$ is not contained in the ideal of $C$ generated by $\RelEnts$
and $\trace(\bx\by) - v$.  Therefore, for these choices of $k$, $R$
has infinitely many equivalence classes of semisimple $3$-dimensional
representations. (All of the calculations discussed in this section
were performed on computers with $\leq 8$ GB RAM.)
\end{note}

\begin{note}
For a second example, consider the case when
\[ R \; = \; k\{ X , Y \} \big/ \langle f_1, \, f_2 \rangle, \]
for
\[ f_1 \, = \, (XY - YX)X - X(XY - YX)-2X, \; f_2 \, = \, (XY - YX)Y - Y(XY -
YX)+2Y. \]
It is well known that $R$, now, is isomorphic to the enveloping algebra of
$\sltwo(k)$. It is also well known, when $k$ has characteristic zero, that
$\sltwo(k')$ has exactly one irreducible representation (up to equivalence) in
each finite dimension.

We used Macaulay2 to implement the algorithm for $k = \Q$, $\F_2$, $\F_3$,
$\F_5$, and $\F_7$. The procedure showed that there are only finitely many
equivalence classes of $3$-dimensional representations when $k = \Q$, $\F_5$,
and $\F_7$. For $k = \F_2$, the procedure found that $1$ is not in the ideal
of $C$ generated by $\RelEnts$ and $\trace(\bx)-v$, thus showing the existence
of infinitely many distinct equivalence classes of semisimple $3$-dimensional
representations. For $k = \F_3$, the procedure found that $1$ is not contained
in the ideal of $C$ generated by $\RelEnts$ and
\[ (\text{the degree-$2$ coefficient of the characteristic polynomial of
$\bx\by$}) - v ,\]
demonstrating the existence of infinitely many distinct equivalence classes of
semisimple $3$-dimensional representations.
\end{note}

\end{document}